\documentclass[12pt,reqno]{amsart}

\usepackage{graphicx}
\graphicspath{ {./images/} }
\usepackage{bm}

\usepackage{url}

\usepackage{amssymb}

\usepackage{amsfonts}

\usepackage{amsmath}

\usepackage[all]{xy}

\newtheorem{theorem}{Theorem}[section]

\newtheorem{corollary}[theorem]{Corollary}

\newtheorem{lemma}[theorem]{Lemma}

\newtheorem{Definition}[theorem]{Definition}

\newtheorem{Example}[theorem]{Example}

\newtheorem{Remark}[theorem]{Remark}

\newenvironment{remark}{\begin{Remark}\begin{em}}{\end{em}\end{Remark}}
\newenvironment{example}{\begin{Example}\begin{em}}{\end{em}\end{Example}}

\newcommand{\T}{\mathrm{UT}} % positive-diagonal upper triangular group

\newcommand{\CP}{\mathrm{CP}}

\newcommand{\UP}{\mathrm{UP}}

\DeclareMathOperator{\diag}{diag}

\address{Huajun Huang \\ Department of Mathematics and Statistics, Auburn University, Auburn, AL 36849, USA}
\email{huanghu@auburn.edu}

\address{Sejong Kim \\ Department of Mathematics, Chungbuk National University, Cheongju 28644, Korea}
\email{skim@chungbuk.ac.kr}

\address{Vatsalkumar N. Mer \\ Institute for Industrial and Applied Mathematics, Chungbuk National University, Cheongju 28644, Korea}
\email{vnm232657@gmail.com}

\begin{document}

\author[Huajun Huang, Sejong Kim and Vatsalkumar N. Mer]{Huajun Huang, Sejong Kim and Vatsalkumar N. Mer}

\title[Power Mean, Karcher Mean, and Log-Euclidean Mean on $\CP_n$]{Power Mean, Karcher Mean, and Log-Euclidean Mean on Positive Constant Upper Triangular Matrices}

\date{}

\begingroup
\begin{abstract}

We study weighted power, Karcher, and log-Euclidean means on the Lie
group $\CP_n$ of real upper triangular matrices with a constant positive
diagonal. Building on the general triangular-group theory, we use the
central decomposition $\CP_n\cong\mathbb R_+\times\UP_n$ to separate
the scalar power mean from a unipotent mean with reweighted coefficients.
We prove that the power mean extends real-analytically through order
zero, where it equals the Karcher mean. We establish exterior-power
covariance for the Karcher and log-Euclidean means, and tensor-product
identities for all three means, distinguishing product weights from
paired inputs. Finite nilpotent functional calculus yields concise
formulas on $\CP_2$ and $\CP_3$, where the Karcher and log-Euclidean
means coincide. In dimension four, their difference is exactly a central
cubic commutator, giving a necessary and sufficient equality criterion
and an explicit counterexample. In every dimension, the same commutator
is the leading term of their small-dispersion difference.
\end{abstract}
\endgroup
\maketitle

\vspace{5mm}



\begingroup
\section{Introduction}

Matrix and operator means provide nonlinear averages on the cone
$\mathbb P_n$ of positive definite matrices and on cones of positive
invertible operators. For $A_1,\ldots,A_m\in\mathbb P_n$ and a
probability vector $\omega=(w_1,\ldots,w_m)$, the weighted Karcher mean
is the unique minimizer
\[
 \Lambda(\omega;A_1,\ldots,A_m)
 =\underset{X\in\mathbb P_n}{\arg\min}
   \sum_{i=1}^m w_i d_R(X,A_i)^2,
 \qquad
 d_R(A,B)=\|\log(A^{-1/2}BA^{-1/2})\|_2,
\]
where $\|\cdot\|_2$ denotes the Hilbert--Schmidt norm. Equivalently,
it solves the Karcher equation \cite{LL14}. The weighted power means
introduced by Lim and P\'alfia \cite{LimPalfia} solve
\begin{equation}\label{E:power}
 X=\sum_{i=1}^m w_i(X\sharp_t A_i),\qquad 0<t\leq1,
\end{equation}
where $A\sharp_t B=A^{1/2}(A^{-1/2}BA^{-1/2})^tA^{1/2}$.
On $\mathbb P_n$, these means converge to the Karcher mean as
$t\to0^+$, with additional L\"owner-order properties
\cite{LimPalfia,LL14}.

We consider the analogous algebraic means on triangular Lie groups.
Let $\T_n$ be the group of real upper triangular matrices with positive
diagonal, let $\UP_n$ be its unipotent subgroup, and let $\CP_n$ be the
subgroup with a constant positive diagonal. Every $A\in\CP_n$ has the
unique central factorization
\[
 A=a\widehat A,\qquad a>0,\quad\widehat A\in\UP_n,
 \qquad \CP_n\cong\mathbb R_+\times\UP_n.
\]
This decomposition separates scalar averaging from nilpotent matrix
operations. For the power mean the unipotent factor uses
$t$-dependent weights, whereas the Karcher and log-Euclidean means
retain the original weights and use the geometric average of the scalar
factors.

Choi, Kim, and Lim obtained finite binomial formulas for geometric
means of unipotent matrices and studied their low-dimensional relation
to the log-Euclidean mean \cite{CKL}. Lawson developed weighted Karcher
means on unipotent and related Lie groups \cite{Lawson24,Lawson25},
and operator means on $\CP_n$ were studied in \cite{HKM}.
Huang and Lawson  established existence, uniqueness,
analyticity, and recursive formulas for Karcher means on $\T_n$
\cite{HuangLawsonKarcher2026}. Their treatment of weighted power means
on $\T_n$ includes the corresponding structural properties, convergence
to the Karcher mean, and Lie--Trotter formulas
\cite{HuangLawsonPower2026}.
We use these general results as background and cite them without
reproducing their proofs. The emphasis here is on the additional
descriptions supplied by the central splitting and the nilpotent
structure of $\CP_n$.

The results are organized around four features of this setting.
\begin{enumerate}
\item \emph{Scalar splitting and the order parameter.}
The power mean has scalar factor $\alpha_t$ given in
\eqref{power alpha}, and its unipotent factor is
computed with weights $q_i(t)=w_i a_i^t/\sum_jw_j a_j^t$.
Theorem~\ref{T:power to Karcher} strengthens the known convergence
at order zero to real-analytic continuation through $t=0$.

\item \emph{Exterior powers and tensor products.}
Theorem~\ref{thm:exterior-means} gives exterior-power covariance for
$\Lambda$ and $\mathrm{LE}$ and a scalar counterexample for $P_t$.
Theorem~\ref{thm:CP-tensor} gives tensor-product multiplicativity of
all three means for product weights. For paired inputs,
$\Lambda$ and $\mathrm{LE}$ remain multiplicative, while $P_t$
generally does not (Theorem~\ref{thm:paired-tensor}). These are
structural consequences of the defining equations and the compatible
exponential and logarithm maps.

\item \emph{Closed formulas in dimensions two and three.}
Theorems~\ref{T:power-CP3} and~\ref{T:LE-CP3} express the three means
using normalized entries and a single weighted covariance term.
The equality $\Lambda=\mathrm{LE}$ on $\CP_3$ yields the exact
common-power identity in Corollary~\ref{P:Karcher} and the explicit
power-mean correction in Corollary~\ref{cor:CP3-power-difference}.

\item \emph{The first obstruction to equality.}
For $A_i=a_i e^{Y_i}\in\CP_4$, set
$\alpha=\prod_i a_i^{w_i}$ and $M=\sum_iw_iY_i$. Then
Theorem~\ref{thm:CP4-KLE-obstruction} proves the exact identity
\[
 \Lambda(\omega;\mathbf A)-\mathrm{LE}(\omega;\mathbf A)
 =\frac{\alpha}{12}\sum_iw_i[Y_i,[M,Y_i]].
\]
The right-hand side is central and gives a necessary and sufficient
condition for equality. Example~\ref{ex:CP4-strict} shows that four
is the first dimension where the means can differ.
Corollary~\ref{cor:CPn-small-dispersion} identifies the same cubic
term in every dimension.
\end{enumerate}

Section~2 recalls the geometric, Karcher, and log-Euclidean means.
Section~3 separates the known power-mean theory from analytic
continuation, exterior powers, and tensor products. Section~4 gives
the low-dimensional formulas and their consequences. Section~5
establishes the dimension-four obstruction and the general cubic
expansion.
\endgroup

\section{{Geometric, Karcher, and log-Euclidean means}}

\subsection{Functions on upper triangular Lie groups}

A unipotent matrix is an upper triangular matrix whose main diagonal entries are all $1$.
Let $\operatorname{UP}_n$ denote the Lie group of $n \times n$ unipotent matrices over $\mathbb{R}$.
Its Lie algebra, denoted by $\mathfrak{up}_n$, consists of all $n \times n$ nilpotent
upper triangular matrices over $\mathbb{R}$.

Let $\operatorname{CP}_n$ denote the Lie group of $n \times n$ upper triangular matrices over $\mathbb{R}$ whose diagonal entries are all equal to the same positive constant.
Its Lie algebra $\mathfrak{cp}_n$ consists of upper triangular matrices whose diagonal
entries are all equal. More explicitly,
\begin{align*}
\operatorname{CP}_n	&= \mathbb{R}_{+} I_{n} \cdot \operatorname{UP}_n
\cong \mathbb{R}_{+} \times \operatorname{UP}_n, \\
\mathfrak{cp}_n	&= \mathbb{R} I + \mathfrak{up}_n,
\end{align*}
where $I_{n} \in \operatorname{UP}_n$ denotes the $n \times n$ identity matrix.

More generally, let $\T_n$ denote the Lie group of all
$n\times n$ upper triangular matrices over $\mathbb{R}$ with positive
diagonal entries, which are not necessarily equal.
Let $\mathfrak{t}_n$ denote the Lie algebra of $\T_n$
consisting of all $n \times n$ upper triangular matrices over $\mathbb{R}$.
Then
\[
\operatorname{UP}_n
\subset
\operatorname{CP}_n
\subset
\T_n,
\qquad
\mathfrak{up}_n
\subset
\mathfrak{cp}_n
\subset
\mathfrak{t}_n.
\]

The exponential map
\begin{equation}\label{E:exp}
	\exp : \mathfrak{t}_n \longrightarrow \T_n,
	\qquad
	\exp X
	= \sum_{k=0}^{\infty} \frac{X^k}{k!},
	\quad
	X \in \mathfrak{t}_n,
\end{equation}
is a bijection whose inverse is the logarithm map
\[
\log : \T_n \longrightarrow \mathfrak{t}_n.
\]
Moreover, these maps restrict to mutually inverse bijections
\begin{displaymath}
\begin{split}
\exp : \mathfrak{cp}_n \longrightarrow \operatorname{CP}_n,
\qquad &
\log : \operatorname{CP}_n \longrightarrow \mathfrak{cp}_n, \\
\exp : \mathfrak{up}_n \longrightarrow \operatorname{UP}_n,
\qquad &
\log : \operatorname{UP}_n \longrightarrow \mathfrak{up}_n.
\end{split}
\end{displaymath}
One can naturally define the $t$-power of $A \in \T_n$ for $t \in \mathbb{R}$ by
\begin{displaymath}
A^{t} = \exp (t \log A).
\end{displaymath}
This map plays a very important role in this paper.

\subsection{The weighted geometric mean and the Karcher mean on $\CP_n$}

The weighted geometric mean of $A, B \in \T_n$ is given by
\begin{equation} \label{E:geomean}
t \in [0,1] \mapsto A \sharp_{t} B = A^{1/2} (A^{-1/2} B A^{-1/2})^{t} A^{1/2} = A  (A^{-1} B)^{t}.
\end{equation}
When $A$ and $B$ are elements of $\CP_n$ (resp. $\UP_n$), the matrix
$A \sharp_{t} B$ is also in $\CP_n$ (resp. $\UP_n$).

\begingroup
The following algebraic properties are consequences of analytic functional calculus and are also discussed on $\CP_n$ in \cite{HKM}.
\endgroup

\begin{theorem} \label{T:t-geomean}
Let $A, B, U, V \in \T_n$ and $s, t, u \in \mathbb{R}$.
\begin{itemize}
\item[(1)] $A \sharp_t B = A (A^{-1} B)^t = (A B^{-1})^{1-t} B$.
\item[(2)] $A \sharp_{t} B = A^{1-t} B^{t}$ if $A$ and $B$ commute.
\item[(3)] $(\alpha A) \sharp_{t} (\beta B) = \alpha^{1-t} \beta^{t} (A \sharp_{t} B)$ for any $\alpha, \beta > 0$.
\item[(4)] $(A \sharp_{t} B )^{-1} = A^{-1} \sharp_{t} B^{-1}$.
\item[(5)] $A \sharp_{t} B = B \sharp_{1-t} A$.
\item[(6)] $(U A V) \sharp_{t} (U B V) = U (A \sharp_{t} B) V$.
\item[(7)] $(A \sharp_{s} B) \sharp_{u} (A \sharp_{t} B) = A \sharp_{(1-u) s + u t} B$.
\end{itemize}
\end{theorem}

In the remaining of this section, we focus on the Lie group $\CP_n$.

We compute $A \sharp_t B$ for $A, B \in \CP_n$ explicitly.
Let $A = a I_{n} + N \in \CP_n$ where $a \in \mathbb{R}_+$ and $N$ is a nilpotent upper triangular matrix.
The binomial expansion of $A^{t}$ is
\begin{equation} \label{E:A^t-binomial}
A^{t} = (a I_{n} + N)^{t} = \sum_{i=0}^{n-1} a^{t-i} \binom{t}{i} N^{i},
\end{equation}
where
$$ \binom{t}{i} = \frac{t(t-1)(t-2)\cdots (t-i+1)}{i!}, \qquad
\binom{t}{0} := 1. $$

\begingroup
For $A=a\widehat A$ and $B=b\widehat B$ in $\CP_n$,
with $\widehat A,\widehat B\in\UP_n$, the nilpotent matrix
$N=\widehat A^{-1}\widehat B-I_n$ satisfies $N^n=0$. Hence
\begin{equation}\label{E:bi-expansion}
 A\sharp_tB=a^{1-t}b^t\widehat A
 \sum_{k=0}^{n-1}\binom tk
 (\widehat A^{-1}\widehat B-I_n)^k,
 \qquad t\in\mathbb R.
\end{equation}
This is an exact \emph{finite} identity with no restriction on the
scalar ratio $b/a$. It is the central-scaled version of the unipotent
binomial expansion \cite{CKL}; see also \cite{HKM}.

\endgroup

We now turn to the multivariate Karcher mean on $\CP_n$.

\begingroup
Throughout, write
\[
 \Delta_m=\{\omega\in\mathbb R^m:w_i\geq0,\ \sum_iw_i=1\},\qquad
 \Delta_m^{\circ}=\{\omega\in\Delta_m:w_i>0\text{ for every }i\}.
\]
All means are defined for $\omega\in\Delta_m$ by deleting zero-weight
inputs. Analytic dependence on weights is understood on
$\Delta_m^{\circ}$.
\endgroup

From \cite{Lawson24} the Karcher (exponential) mean $\Lambda(\omega; A_{1}, \dots, A_{m})$ of $A_{1}, \dots, A_{m} \in \UP_{n}$ is a unique solution $X \in \UP_{n}$ to the (Karcher) equation
\begin{equation} \label{E:Karcher-eq}
\sum_{i=1}^{m} w_{i} \log (X^{-1} A_{i}) = 0.
\end{equation}
Furthermore, it has been shown in \cite{Lawson25} that the sequence $\{ X_{k} \}_{k \geq 1}$ given by
\begin{displaymath}
X_{k+1} = X_{k} \exp \left[ \sum_{i=1}^{m} w_{i} \log (X_{k}^{-1} A_{i}) \right]
\end{displaymath}
with any $X_{1} \in \UP_{n}$ converges to $\Lambda(\omega; A_{1}, \dots, A_{m})$ in at most $n-1$ steps.
By joint homogeneity of the Karcher mean, we can define the Karcher mean of $A_{1}, \dots, A_{m} \in \CP_{n}$:
\begin{equation} \label{E:K-mean}
\Lambda(\omega; A_{1}, \dots, A_{m}) = a_{1}^{w_{1}} \cdots a_{m}^{w_{m}} \Lambda(\omega; \widehat{A}_{1}, \dots, \widehat{A}_{m})
\end{equation}
where $A_{1} = a_{1} \widehat{A}_{1}, \dots, A_{m} = a_{m} \widehat{A}_{m}$ for some positive constants $a_{i}$ and $\widehat{A}_{1}, \dots, \widehat{A}_{m} \in \UP_{n}$.

\begin{remark}
For $A_{1}, \dots, A_{m} \in \CP_{n}$, the Karcher mean $\Lambda(\omega; A_{1}, \dots, A_{m})$ given by \eqref{E:K-mean} also satisfies the Karcher equation \eqref{E:Karcher-eq}.
Indeed, let $X = \Lambda(\omega; A_{1}, \dots, A_{m})$ for $A_{1} = a_{1} \widehat{A}_{1}, \dots, A_{m} = a_{m} \widehat{A}_{m}$ for some positive constants $a_{i}$ and $\widehat{A}_{1}, \dots, \widehat{A}_{m} \in \UP_{n}$. Then $X = \alpha \widehat{X}$, where
\begin{displaymath}
\widehat{X} = \Lambda(\omega; \widehat{A}_{1}, \dots, \widehat{A}_{m}) \quad \textrm{and} \quad \alpha = \prod_{i=1}^{m} a_{i}^{w_{i}}.
\end{displaymath}
Also, $\log (X^{-1} A_{i}) = \log (\alpha^{-1} a_{i}) I_{n} + \log (\widehat{X}^{-1} \widehat{A}_{i})$ for each $i$.
Thus,
\begin{displaymath}
\sum_{i=1}^{m} w_{i} \log (X^{-1} A_{i}) = \sum_{i=1}^{m} w_{i} \log (\alpha^{-1} a_{i}) I_{n} + \sum_{i=1}^{m} w_{i} \log (\widehat{X}^{-1} \widehat{A}_{i}) = 0.
\end{displaymath}
\end{remark}

\begingroup
\begin{lemma}[\cite{HuangLawsonKarcher2026}]
Let $t \in [0,1]$. For $m=2$ and $\omega = (1-t,t)$, the Karcher mean coincides with the weighted geometric mean: for $A_{1}, A_{2} \in \CP_n$
\[
\Lambda (\omega; A_1, A_2) = A_1 \sharp_t A_2.
\]
\end{lemma}

\endgroup

For a square matrix $A$, let $\diag(A)$ denote the diagonal matrix formed by the diagonal of $A$.
\begingroup
Restriction to a contiguous principal block is an algebra homomorphism
on upper triangular matrices. Thus, for a one-variable function $f$
analytic near the spectrum of $A$, one has
\begin{align}
 \diag(f(A))&=f(\diag A),\label{T diag fun}\\
 f(A)[i:j]&=f(A[i:j]),\qquad 1\leq i\leq j\leq n.
 \label{T_block_fun}
\end{align}
Here $A[i:j]$ is the principal block on the consecutive indices
$i,\ldots,j$; see \cite{HuangLawsonKarcher2026,HuangLawsonPower2026}.
Applying these identities to the Karcher equation gives the following
known diagonal and block properties.
\endgroup

\begingroup
\begin{lemma}[\cite{HuangLawsonKarcher2026}]\label{thm: mean property}
Let \( \omega = (w_1, \dots, w_m) \in \Delta_{m} \) and \( A_1, \dots, A_m \in \CP_n \).
Then the following properties for $X = \Lambda(\omega; A_1, \dots, A_m)$ hold:
\begin{enumerate}
\item $\diag(X) = \diag (A_1^{w_1} \cdots A_m^{w_m})$.

\item For $1 \le i \le j \le n$, the principal submatrix $X[i:j]$ is an $\omega$-weighted Karcher mean of	$A_1[i:j], \ldots, A_m[i:j]$.
	In particular, the $(i,j)$-entry of $X$ depends only on $\omega$ and the principal submatrices $A_1[i:j], \ldots, A_m[i:j]$.
	\end{enumerate}
\end{lemma}

\endgroup

\begin{theorem} \label{T:Karcher-properties}
Let \( \omega = (w_1, \dots, w_m) \in \Delta_{m} \) and \( \mathbb{A} = (A_1, \dots, A_m) \in (\CP_n)^{m} \).
\begin{itemize}
  \item[(1)] $\Lambda(\omega; \mathbb{A}) = A_{1}^{w_{1}} \cdots A_{m}^{w_{m}}$ if the $A_{i}$'s commute;
  \item[(2)] $\Lambda(\omega; c_{1} A_1, \dots, c_{m} A_m) = c_{1}^{w_{1}} \cdots c_{m}^{w_{m}} \Lambda(\omega; \mathbb{A})$ for any $c_{1}, \dots, c_{m} > 0$;
  \item[(3)] $\Lambda(\omega_{\sigma}; \mathbb{A}_{\sigma}) = \Lambda(\omega; \mathbb{A})$ for any permutation $\sigma$ of $\{ 1, \cdots , m \}$, where
  \begin{displaymath}
  \begin{split}
  \omega_{\sigma} & := (w_{\sigma(1)}, \dots, w_{\sigma(m)}) \in \Delta_{m}, \\
  \mathbb{A}_{\sigma} & := (A_{\sigma(1)}, \dots, A_{\sigma(m)}) \in (\CP_{n})^{m};
  \end{split}
  \end{displaymath}
  \item[(4)] $\Lambda (\omega; U A_1 V, \ldots, U A_m V) = U \Lambda (\omega; \mathbb{A}) V$ for any $U, V \in \CP_n$;
  \item[(5)] $\Lambda(\omega; \mathbb{A}^{-1}) = \Lambda(\omega; \mathbb{A})^{-1}$.

\end{itemize}
\end{theorem}

\begingroup
\begin{proof}
The two-sided covariance in (4) is established on $\T_n$ in
\cite{HuangLawsonKarcher2026}, so we omit its proof. Parts (1)--(3)
follow from the defining equation, the central splitting, and uniqueness.
For (5), if $X=\Lambda(\omega;\mathbf A)$, then
\[
 \sum_iw_i\log(XA_i^{-1})
 =-X\left(\sum_iw_i\log(X^{-1}A_i)\right)X^{-1}=0.
\]
Uniqueness identifies $X^{-1}$ as the mean of the inverse inputs.
\end{proof}
\endgroup

\subsection{The log-Euclidean mean}

\begingroup
For $A_1, \dots, A_m \in \T_n$ and $\omega = (w_1, \dots, w_m) \in \Delta_m$, we define the log-Euclidean mean by
\endgroup

\begin{equation} \label{E:log-Euclidean}
\mathrm{LE}(\omega; A_1, \dots, A_m)
= \exp \left( \sum_{k=1}^m w_k \log A_k \right).
\end{equation}
Suppose now that $A_1, \dots, A_m \in \CP_n$.
Write $A_k = a_k \widehat A_k$, where
$a_k > 0$, and $\widehat A_k \in \UP_n$ for $k = 1, \dots, m$.
Since $a_k I_{n}$ commutes with $\widehat A_k$, we have
$\log A_k = (\log a_k) I_{n} + \log \widehat A_k$.
Consequently,
\begin{align}
\mathrm{LE}(\omega; A_1, \dots, A_m)
& = \exp \left[ \left( \sum_{k=1}^m w_k \log a_k \right) I_{n}
+ \sum_{k=1}^m w_k \log \widehat A_k \right] \notag \\
& = \left( \prod_{k=1}^m a_k^{w_k} \right)
\mathrm{LE}(\omega; \widehat A_1, \dots, \widehat A_m). \label{E:LE-CP}
\end{align}

\begin{theorem} \label{T:LE-properties}
Let $\omega = (w_1, \dots, w_m) \in \Delta_m$ and $\mathbb{A} = (A_1, \dots, A_m) \in (\T_n)^m$. Then the log-Euclidean mean satisfies the following properties.
\begin{enumerate}
 \item[(1)] $\mathrm{LE}(\omega; \mathbb{A}) = A_1^{w_1} \cdots A_m^{w_m}$ if the $A_i$'s commute.
 \item[(2)] $\mathrm{LE}(\omega; a_1 A_1, \dots, a_m A_m) = a_1^{w_1} \cdots a_m^{w_m}	\mathrm{LE}(\omega; \mathbb{A})$ for any $a_1, \dots, a_m > 0$.
 \item[(3)] $\mathrm{LE}(\omega_\sigma; \mathbb{A}_{\sigma}) = \mathrm{LE}(\omega; \mathbb{A})$ for every permutation $\sigma$ on $\{ 1, \dots, m \}$.
 \item[(4)] $\mathrm{LE} \bigl( \omega; P^{-1} A_1 P, \dots, P^{-1} A_m P \bigr)
	= P^{-1} \mathrm{LE}(\omega; \mathbb{A}) P$ for every $P \in \T_n$.
 \item[(5)] $\mathrm{LE}(\omega; \mathbb{A}^{-1}) = \mathrm{LE}(\omega; \mathbb{A})^{-1}$.
 \item[(6)]	$\mathrm{LE}(\omega; A_1^r, \dots, A_m^r) = \mathrm{LE}(\omega; \mathbb{A})^r$ for every $r \in \mathbb{R}$.
 \item[(7)]	$\mathrm{LE}(\omega; \mathbb{A})[i:j] = \mathrm{LE} \bigl( \omega; A_1[i:j], \dots, A_m[i:j] \bigr)$ for $1 \leq i \leq j \leq n$.
 In particular, $\displaystyle \diag \bigl( \mathrm{LE}(\omega; \mathbb{A}) \bigr) = \prod_{k=1}^m \diag(A_k)^{w_k}$.
\end{enumerate}
\end{theorem}

\begin{proof}
Properties (1)-(3) are obvious from \eqref{E:log-Euclidean} and \eqref{E:LE-CP}.
Using $\log(A_k^{-1}) = - \log A_k$ and $\log (A_k^r) = r \log A_k$ for any $r \in \mathbb{R}$, we obtain (5) and (6).
\begin{enumerate}		
\item[(4)] The analytic functional calculus gives $\log(P^{-1} A_k P) =	P^{-1} (\log A_k) P$. Thus,
\begin{align*}
\mathrm{LE} \bigl( \omega; P^{-1} A_1 P, \dots, P^{-1} A_m P \bigr)
& = \exp \left(	\sum_{k=1}^m w_k P^{-1} (\log A_k) P \right) \\
& =	P^{-1} \exp \left( \sum_{k=1}^m w_k \log A_k \right) P = P^{-1} \mathrm{LE}(\omega; \mathbb{A}) P.
\end{align*}
		
\item[(7)] For an upper triangular matrix $A$ and an analytic one-variable function $f$, one has $f(A)[i:j] = f(A[i:j])$. Consequently,
\[
(\log A_k)[i:j]	= \log(A_k[i:j]).
\]
	Applying this equality first to the logarithm and then to the exponential gives
\begin{align*}
\mathrm{LE}(\omega; \mathbb{A})[i:j]
& = \exp \left(	\sum_{k=1}^m w_k \log A_k \right)[i:j] \\
& =	\exp \left(	\sum_{k=1}^m w_k \log(A_k[i:j])	\right) \\
& =	\mathrm{LE}	\bigl( \omega; A_1[i:j], \dots, A_m[i:j] \bigr).
\end{align*}
Taking $i=j$ yields the last assertion.
\end{enumerate}
\end{proof}

\begingroup
The log-Euclidean mean also occurs in Lie--Trotter formulas. The
positive definite theory is treated in \cite{HK17}. For triangular
groups, the all-dimensional Karcher and power-mean formulas are known
from \cite{HuangLawsonPower2026} and are recorded together in
Theorem~\ref{T:power-Lie-Trotter} below.
\endgroup

\begingroup
\section{Power means: scalar splitting, limits, and structural identities}

\subsection{Known theory and the scalar--unipotent splitting}

For $\mathbf A=(A_1,\ldots,A_m)\in(\CP_n)^m$,
$\omega\in\Delta_m$, and $0<t\leq1$, the power-mean equation is
\begin{equation}\label{E:p-mean}
 X=\sum_{i=1}^m w_i(X\sharp_t A_i).
\end{equation}
Using $X\sharp_t A_i=X(X^{-1}A_i)^t$, this is equivalent to
\begin{equation}\label{E:p-mean-2}
 \sum_{i=1}^m w_i(X^{-1}A_i)^t=I_n.
\end{equation}
The following is the restriction of the existence and analyticity
theorem on $\T_n$ in \cite{HuangLawsonPower2026}. Its diagonal formula
ensures that the solution belongs to $\CP_n$ when the inputs do.

\begin{theorem}[\cite{HuangLawsonPower2026}]\label{T:power-mean}
For $0<t\leq1$ and $\mathbf A\in(\CP_n)^m$, equation
\eqref{E:p-mean} has a unique solution
$P_t(\omega;\mathbf A)\in\CP_n$. For fixed $t$, this solution
depends real-analytically on the input entries and on
$\omega\in\Delta_m^{\circ}$.
\end{theorem}

We also record the structural identities from
\cite{HuangLawsonPower2026}, omitting their proofs.

\begin{theorem}[\cite{HuangLawsonPower2026}]
\label{T:properties of p-mean}
Let $\mathbf A\in(\CP_n)^m$, $\omega\in\Delta_m$, and $0<t\leq1$.
\begin{enumerate}
\item If the inputs commute pairwise, then
\[
 P_t(\omega;\mathbf A)=\left(\sum_iw_iA_i^t\right)^{1/t}.
\]
\item For $c_i>0$, put $s_t=(\sum_iw_ic_i^t)^{1/t}$ and
$v_i=w_ic_i^t/\sum_jw_jc_j^t$. Then
\[
 P_t(\omega;c_1A_1,\ldots,c_mA_m)=s_tP_t(v;\mathbf A).
\]
\item Simultaneous permutation of the inputs and weights leaves
the mean unchanged.
\item For $U,V\in\CP_n$,
\[
 P_t(\omega;UA_1V,\ldots,UA_mV)=UP_t(\omega;\mathbf A)V.
\]
\item If $m\geq2$ and $w_m<1$, then
\[
 P_t(\omega;A_1,\ldots,A_{m-1},X)=X
 \ \Longleftrightarrow\ 
 P_t(\widehat\omega;A_1,\ldots,A_{m-1})=X,
\]
where $\widehat\omega=(w_1,\ldots,w_{m-1})/(1-w_m)$.
\item For $1\leq p\leq q\leq n$,
\[
 P_t(\omega;\mathbf A)[p:q]
 =P_t(\omega;A_1[p:q],\ldots,A_m[p:q]).
\]
If $a_i$ is the constant diagonal entry of $A_i$, the mean has
constant diagonal
\begin{equation}\label{power alpha}
 \alpha_t=\left(\sum_iw_i a_i^t\right)^{1/t}.
\end{equation}
\end{enumerate}
For boundary weights, these statements follow by deleting the
zero-weight inputs.
\end{theorem}

Write $A_i=a_i\widehat A_i$ with $\widehat A_i\in\UP_n$ and set
\begin{equation}\label{power adjusted w}
 q_i(t)=\frac{w_i a_i^t}{\sum_jw_j a_j^t},\qquad
 q(t)=(q_1(t),\ldots,q_m(t)).
\end{equation}
The homogeneity identity above gives the exact factorization
\begin{equation}\label{eq:CP-splitting-limit}
 P_t(\omega;\mathbf A)
 =\alpha_tP_t(q(t);\widehat A_1,\ldots,\widehat A_m).
\end{equation}
Thus computation on $\CP_n$ reduces to scalar averaging and a
unipotent equation. Notice that the unipotent weights depend on $t$
unless the scalar factors coincide on the support of $\omega$.

\subsection{Analytic continuation at zero and Lie--Trotter limits}

Convergence of $P_t$ to the Karcher mean is known on the larger
group $\T_n$ \cite{HuangLawsonPower2026}. The constant-diagonal
splitting yields the following stronger parameter statement.

\begin{theorem}[Analytic continuation through order zero]
\label{T:power to Karcher}
For $\mathbf A\in(\CP_n)^m$ and $\omega\in\Delta_m^{\circ}$,
the map $t\mapsto P_t(\omega;\mathbf A)$, initially defined for
$0<t\leq1$, extends real-analytically to a neighborhood of zero.
Its value at zero is $\Lambda(\omega;\mathbf A)$. In particular,
\begin{equation}\label{eq:CP-power-limit}
 \lim_{t\to0^+}P_t(\omega;\mathbf A)=\Lambda(\omega;\mathbf A).
\end{equation}
\end{theorem}

\begin{proof}
The functions $q_i(t)$ in \eqref{power adjusted w} are analytic and
$q_i(0)=w_i$. Also, $\alpha_t$ extends analytically to zero because
$\log(\sum_iw_i e^{t\log a_i})$ vanishes there, and
$\alpha_0=\prod_i a_i^{w_i}$. By \eqref{eq:CP-splitting-limit}, it
suffices to treat fixed unipotent inputs $U_i$ with analytic positive
weights $q_i(t)$ summing to one.

Proceed by induction on $n$, starting with $n=1$. Principal-block
consistency and the induction hypothesis determine analytic leading
and trailing blocks of size $n-1$; they agree on their overlap by
uniqueness. Assemble them into $Y(t)\in\UP_n$ with $Y(t)_{1n}=0$,
and put $B_i(t)=Y(t)^{-1}U_i$. Since $E_{1n}$ annihilates strictly
upper triangular matrices on both sides, finite binomial calculus gives
\[
 ((Y(t)+xE_{1n})^{-1}U_i)^t=B_i(t)^t-txE_{1n}.
\]
For $t>0$, the remaining entry is therefore
\begin{equation}\label{eq:CP-entry-analytic}
 x(t)=\sum_iq_i(t)
       \left[\frac{B_i(t)^t-I_n}{t}\right]_{1n}.
\end{equation}
The quotient extends jointly analytically to $t=0$, with value
$\log B_i(0)$, since $B^t=\exp(t\log B)$. Hence $x(t)$ and the
entire solution extend analytically through zero. At zero the same
entry equation becomes
$x(0)=\sum_iw_i[\log(Y(0)^{-1}U_i)]_{1n}$, so
$Y(0)+x(0)E_{1n}$ solves the Karcher equation. Uniqueness
\cite{Lawson24,HuangLawsonKarcher2026} identifies this limit.
Finally use \eqref{E:K-mean} to restore the scalar factors.
\end{proof}

For completeness, we state the Lie--Trotter formulas already proved
on $\T_n$ in \cite{HuangLawsonPower2026}; their proofs are omitted.

\begin{theorem}[\cite{HuangLawsonPower2026}]
\label{T:power-Lie-Trotter}
Let $\mathbf A\in(\CP_n)^m$, $\omega\in\Delta_m$, and $0<t\leq1$.
Then
\begin{equation}\label{eq:CP-general-LT}
 \lim_{p\to0}P_t(\omega;A_1^p,\ldots,A_m^p)^{1/p}
 =\mathrm{LE}(\omega;\mathbf A)
 =\lim_{p\to0}\Lambda(\omega;A_1^p,\ldots,A_m^p)^{1/p}.
\end{equation}
\end{theorem}

\subsection{Exterior powers}

For $1\leq k\leq n$, let $C_k(A)=\bigwedge^k A$ be the
$k$th multiplicative compound matrix in the lexicographically ordered
exterior basis, and set $N=\binom nk$. By Cauchy--Binet,
$C_k(AB)=C_k(A)C_k(B)$. In particular, if $A=aU\in\CP_n$, then
$C_k(A)=a^kC_k(U)\in\CP_N$ and $C_k(U)\in\UP_N$.
The differential $\rho_k$ of $C_k$ is linear and satisfies
\begin{equation}\label{eq:compound-log-exp}
 C_k(e^H)=e^{\rho_k(H)},\qquad
 \log C_k(A)=\rho_k(\log A).
\end{equation}
These identities follow by differentiating the exterior action and
using the unique triangular logarithm. The map $\rho_k$ preserves
commutators but is generally not multiplicative.

\begin{theorem}[Exterior-power covariance]\label{thm:exterior-means}
For $\mathbf A\in(\CP_n)^m$, $\omega\in\Delta_m$, and
$C_k(\mathbf A)=(C_k(A_1),\ldots,C_k(A_m))$, one has
\begin{align}
 C_k(\mathrm{LE}(\omega;\mathbf A))
 &=\mathrm{LE}(\omega;C_k(\mathbf A)),\label{eq:exterior-LE}\\
 C_k(\Lambda(\omega;\mathbf A))
 &=\Lambda(\omega;C_k(\mathbf A)).\label{eq:exterior-Karcher}
\end{align}
For $2\leq k\leq n$, the corresponding identity for $P_t$
fails in general, even for scalar matrices and $t=1$.
\end{theorem}

\begin{proof}
Linearity of $\rho_k$ in \eqref{eq:compound-log-exp} proves
\eqref{eq:exterior-LE}. If $X=\Lambda(\omega;\mathbf A)$, applying
$\rho_k$ to its Karcher equation gives
\[
 \sum_iw_i\log(C_k(X)^{-1}C_k(A_i))=0.
\]
Uniqueness proves \eqref{eq:exterior-Karcher}. For the last assertion,
take $A_1=I_n$, $A_2=2I_n$, and $\omega=(1/2,1/2)$.
Since $P_1$ is the arithmetic mean, the two expressions would be
$(3/2)^kI_N$ and $(1+2^k)I_N/2$, which differ for $k\geq2$.
\end{proof}

The positive identities reflect functoriality under Lie-group
representations with compatible logarithms; no exterior-power
additivity is available for the power-mean equation.

\subsection{Tensor products}

Let $\mathbf A=(A_1,\ldots,A_m)\in(\CP_{n_1})^m$ and
$\mathbf B=(B_1,\ldots,B_\ell)\in(\CP_{n_2})^\ell$,
with weights $\omega\in\Delta_m$ and $\mu\in\Delta_\ell$.
Order the pairs $(i,j)$ lexicographically and write
\[
 \mathbf A\otimes\mathbf B=(A_i\otimes B_j)_{i,j},\qquad
 \omega\otimes\mu=(w_i\mu_j)_{i,j}.
\]

\begin{theorem}[Product-weight multiplicativity]\label{thm:CP-tensor}
For $F=\Lambda$, $F=\mathrm{LE}$, or $F=P_t$ with $0<t\leq1$,
\begin{equation}\label{eq:tensor-all}
 F(\omega\otimes\mu;\mathbf A\otimes\mathbf B)
 =F(\omega;\mathbf A)\otimes F(\mu;\mathbf B).
\end{equation}
\end{theorem}

\begin{proof}
For positive-diagonal triangular matrices, uniqueness of the
triangular logarithm gives
\begin{equation}\label{eq:tensor-log}
 \log(A\otimes B)=\log A\otimes I+I\otimes\log B,
 \qquad (A\otimes B)^s=A^s\otimes B^s.
\end{equation}
The two summands in the logarithm commute. This proves the
log-Euclidean identity directly. For the Karcher means $X$ and $Y$
of the separate input families, the weighted logarithm equation for
$X\otimes Y$ is
\[
 \left(\sum_iw_i\log(X^{-1}A_i)\right)\otimes I
 +I\otimes\left(\sum_j\mu_j\log(Y^{-1}B_j)\right)=0.
\]
For the separate power means $X$ and $Y$, the corresponding equation is
\[
 \sum_{i,j}w_i\mu_j((X\otimes Y)^{-1}(A_i\otimes B_j))^t
 =\left(\sum_iw_i(X^{-1}A_i)^t\right)
   \otimes\left(\sum_j\mu_j(Y^{-1}B_j)^t\right)=I.
\]
Uniqueness proves both assertions.
\end{proof}

For equally sized input families, set
$\mathbf A\boxtimes\mathbf B=(A_1\otimes B_1,\ldots,A_m\otimes B_m)$.

\begin{theorem}[Paired tensor products]\label{thm:paired-tensor}
For $\omega\in\Delta_m$ and $F\in\{\Lambda,\mathrm{LE}\}$,
\begin{equation}\label{eq:paired-all}
 F(\omega;\mathbf A\boxtimes\mathbf B)
 =F(\omega;\mathbf A)\otimes F(\omega;\mathbf B).
\end{equation}
The analogous paired identity fails in general for $P_t$.
\end{theorem}

\begin{proof}
Formula \eqref{eq:tensor-log} separates the paired logarithm sum
into its two marginal sums, proving both positive identities exactly
as above. For failure at $t=1$, take
$A_1=B_1=I$, $A_2=B_2=2I$, and $\omega=(1/2,1/2)$, with identities
of the appropriate sizes. The paired arithmetic mean is
$(5/2)I$, while the tensor product of the arithmetic means is
$(9/4)I$.
\end{proof}

\begin{corollary}\label{cor:tensor-limit}
Under the hypotheses of Theorem~\ref{thm:CP-tensor},
\[
 \lim_{t\to0^+}P_t(\omega\otimes\mu;\mathbf A\otimes\mathbf B)
 =\Lambda(\omega;\mathbf A)\otimes\Lambda(\mu;\mathbf B).
\]
For $\mathbf A\in(\CP_n)^m$ and $1\leq k\leq n$, one also has
\[
 \lim_{t\to0^+}C_k(P_t(\omega;\mathbf A))
 =\Lambda(\omega;C_k(\mathbf A)).
\]
\end{corollary}
Both assertions follow immediately from \eqref{eq:CP-power-limit},
continuity, and the preceding structural identities.
\endgroup

\begingroup
\section{Closed formulas in dimensions two and three}

\subsection{Finite functional calculus on $\CP_n$}

General triangular logarithm and exponential formulas, including
confluent divided differences and Parlett recurrences, are available in
\cite{deBoor2005,Higham2008,Parlett1976,HuangLawsonKarcher2026,HuangLawsonPower2026}.
On $\CP_n$ they reduce to the standard finite identities
\begin{align}
 \log(a(I+N))
 &=(\log a)I+\sum_{r=1}^{n-1}\frac{(-1)^{r+1}}{r}N^r,
 \label{eq:CP-finite-log}\\
 \exp(bI+H)
 &=e^b\sum_{r=0}^{n-1}\frac{H^r}{r!},
 \label{eq:CP-finite-exp}
\end{align}
where $a>0$ and $N,H\in\mathfrak{up}_n$.
Together with \eqref{E:A^t-binomial}, these formulas suffice for the
computations below; no separate treatment of repeated diagonal nodes
is required.

For $U=I+uE_{12}+vE_{13}+zE_{23}$ and
$H=hE_{12}+kE_{13}+lE_{23}$,
the relevant specializations are
\begin{equation}\label{eq:CP3-log-exp}
 \log U=\begin{bmatrix}0&u&v-uz/2\\0&0&z\\0&0&0\end{bmatrix},
 \qquad
 e^H=\begin{bmatrix}1&h&k+hl/2\\0&1&l\\0&0&1\end{bmatrix},
\end{equation}
and
\begin{equation}\label{T:CP3-power}
 U^t=\begin{bmatrix}
 1&tu&tv+\tfrac{t(t-1)}2uz\\0&1&tz\\0&0&1
 \end{bmatrix}.
\end{equation}
These are direct finite-series formulas; the corresponding unipotent
binomial calculus is also used in \cite{CKL}.

\subsection{Normalized entries and closed mean formulas}

Write
\begin{equation}\label{T3 ut entries}
 A_i=a_iU_i,\qquad
 U_i=\begin{bmatrix}1&u_i&v_i\\0&1&z_i\\0&0&1\end{bmatrix}
 \in\UP_3,
\end{equation}
so $u_i=(A_i)_{12}/a_i$, $v_i=(A_i)_{13}/a_i$, and
$z_i=(A_i)_{23}/a_i$. For a probability vector $r$, put
\[
 \overline u_r=\sum_i r_i u_i,\quad
 \overline v_r=\sum_i r_i v_i,\quad
 \overline z_r=\sum_i r_i z_i,
\]
and define the scalar weighted covariance
\begin{equation}\label{eq:CP-covariance}
 C_r(u,z)=\sum_i r_i(u_i-\overline u_r)(z_i-\overline z_r)
 =\sum_{i<j}r_ir_j(u_i-u_j)(z_i-z_j).
\end{equation}
This notation isolates the only quadratic correction in dimension three.

\begin{theorem}[Power means]\label{T:power-CP3}
Let $0<t\leq1$, let $\alpha_t$ and $q(t)$ be given by
\eqref{power alpha} and \eqref{power adjusted w}, and abbreviate
$q=q(t)$. For the inputs in \eqref{T3 ut entries},
\begin{equation}\label{power-3x3-form}
 P_t(\omega;\mathbf A)=\alpha_t
 \begin{bmatrix}
 1&\overline u_q&\overline v_q-\tfrac{1-t}{2}C_q(u,z)\\
 0&1&\overline z_q\\0&0&1
 \end{bmatrix}.
\end{equation}
For $A_i=a_i\left[\begin{matrix}1&u_i\\0&1\end{matrix}\right]
\in\CP_2$, the corresponding formula is
\begin{equation}\label{power 2x2 form}
 P_t(\omega;\mathbf A)
 =\alpha_t\begin{bmatrix}1&\overline u_q\\0&1\end{bmatrix}.
\end{equation}
\end{theorem}

\begin{proof}
By \eqref{eq:CP-splitting-limit}, it suffices to solve the
unipotent equation with weights $q$. Writing the unknown mean as
$X=I+xE_{12}+yE_{13}+zE_{23}$,
the first-superdiagonal equations give $x=\overline u_q$ and
$z=\overline z_q$. The $(1,3)$-entry of
$\sum_iq_i(X^{-1}U_i)^t=I$, divided by $t$, is
\[
 0=\sum_iq_i\left(v_i-y-x(z_i-z)
              +\frac{t-1}{2}(u_i-x)(z_i-z)\right).
\]
The linear term in $z_i-z$ sums to zero, and
\eqref{eq:CP-covariance} yields the formula for $y$.
Principal-block consistency gives the $2\times2$ assertion.
\end{proof}

\begin{theorem}[Karcher and log-Euclidean means]\label{T:LE-CP3}
Set
\begin{equation}\label{log-Euclidean alpha}
 \alpha=\prod_i a_i^{w_i}.
\end{equation}
For the inputs in \eqref{T3 ut entries},
\begin{equation}\label{LE T3}
 \Lambda(\omega;\mathbf A)=\mathrm{LE}(\omega;\mathbf A)
 =\alpha\begin{bmatrix}
 1&\overline u_\omega&\overline v_\omega-\tfrac12 C_\omega(u,z)\\
 0&1&\overline z_\omega\\0&0&1
 \end{bmatrix}.
\end{equation}
Equivalently,
\begin{equation}\label{E:LE-CP3}
 \Lambda(\omega;\mathbf A)=\mathrm{LE}(\omega;\mathbf A)
 =\alpha\left[\sum_iw_iU_i
       -\frac12\sum_{i<j}w_iw_j(U_i-U_j)^2\right].
\end{equation}
For $A_i=a_i\left[\begin{matrix}1&u_i\\0&1\end{matrix}\right]
\in\CP_2$, both means equal
$\alpha\left[\begin{matrix}1&\overline u_\omega\\0&1\end{matrix}\right]$.
\end{theorem}

\begin{proof}
The finite logarithm and exponential formulas, followed by
\eqref{E:LE-CP}, give \eqref{LE T3}. Since $(U_i-U_j)^2$ has only
the $(1,3)$-entry $(u_i-u_j)(z_i-z_j)$, this is
\eqref{E:LE-CP3}.

For equality with the Karcher mean, put $Y_i=\log U_i$ and
$Y=\sum_iw_iY_i$. Two-step nilpotence gives
\[
 \log(e^{-Y}e^{Y_i})=-Y+Y_i-\tfrac12[Y,Y_i].
\]
The weighted sum vanishes because $\sum_iw_i[Y,Y_i]=[Y,Y]=0$.
Uniqueness and \eqref{E:K-mean} identify the Karcher mean.
Restriction to a principal $2\times2$ block gives the last assertion.
\end{proof}

For two unipotent inputs, the coincidence with the geometric mean
is established in \cite[Section~5]{CKL}; the formulas above include
arbitrary probability weights and independent positive scalar factors.
They also exhibit the role of two-step nilpotence in the multivariate
Karcher equation.

\subsection{Exact identities and comparisons}

\begin{corollary}[Common powers]\label{P:Karcher}
For $\mathbf A\in(\CP_3)^m$, $\omega\in\Delta_m$, and every
real $p\ne0$,
\begin{equation}\label{Karcher LE}
 \Lambda(\omega;A_1^p,\ldots,A_m^p)^{1/p}
 =\mathrm{LE}(\omega;\mathbf A).
\end{equation}
\end{corollary}
Indeed, apply Theorem~\ref{T:LE-CP3} to the powered inputs and
use the common-power identity of Theorem~\ref{T:LE-properties}.
Thus the Karcher Lie--Trotter formula is exact for every nonzero
$p$ in dimensions two and three. The power-mean Lie--Trotter limit
is already included in Theorem~\ref{T:power-Lie-Trotter}.

\begin{corollary}[Power-mean correction]\label{cor:CP3-power-difference}
If $A_i=aU_i\in\CP_3$ have the same scalar factor $a>0$, then
\begin{equation}\label{eq:CP3-power-difference}
 P_t(\omega;\mathbf A)-\Lambda(\omega;\mathbf A)
 =\frac{at}{2}\sum_{i<j}w_iw_j(U_i-U_j)^2,
 \qquad 0<t\leq1.
\end{equation}
\end{corollary}
Here $q(t)=\omega$, so the identity follows by subtracting
\eqref{E:LE-CP3} from \eqref{power-3x3-form}. The correction has
only a possible $(1,3)$-entry, whose sign may vary with the data.

For comparison, the Wasserstein mean of unipotent inputs in
dimension three satisfies \cite[Theorem~4.3]{KM}
\[
 \Omega(\omega;\mathbf U)
 =\sum_iw_iU_i-\frac14\sum_{i<j}w_iw_j(U_i-U_j)^2.
\]
Consequently,
\[
 \Omega(\omega;\mathbf U)-\Lambda(\omega;\mathbf U)
 =\frac14\sum_{i<j}w_iw_j(U_i-U_j)^2.
\]
These identities compare algebraic corrections for non-Hermitian
matrices and do not assert L\"owner-order inequalities. The next
section determines how the equality $\Lambda=\mathrm{LE}$ first
fails beyond dimension three.
\endgroup

\begingroup
\section{The first difference between Karcher and log-Euclidean mean}
\label{sec:CP-stability-obstruction}

For $Y,Z$ in a Lie algebra write $[Y,Z]=YZ-ZY$.
For $A_k\in\CP_n$, factor
\begin{equation}\label{eq:CP-Y-data}
	A_k=a_k\exp(Y_k),\quad Y_k\in\mathfrak{up}_n,\qquad
	\alpha=\prod_ka_k^{w_k},\quad
	M=\sum_kw_kY_k.
\end{equation}
Thus $\mathrm{LE}(\omega;\mathbf A)=\alpha\exp M$.
Let $\mathfrak n=\mathfrak{up}_4$ and write
$\mathfrak n^2=[\mathfrak n,\mathfrak n]$ and
$\mathfrak n^3=[\mathfrak n,\mathfrak n^2]$.
This Lie algebra is nilpotent of class three, with
$\mathfrak n^3=\mathbb R E_{14}$ and
$[\mathfrak n,\mathfrak n^3]=0$. The class-three
Baker--Campbell--Hausdorff identity is
\begin{equation}\label{eq:CP-BCH3}
	\log(e^Xe^Y)=X+Y+\frac12[X,Y]
	+\frac1{12}[X,[X,Y]]+\frac1{12}[Y,[Y,X]].
\end{equation}
We get the exact dimension-four obstruction. 

\begin{theorem} 
	\label{thm:CP4-KLE-obstruction}
	Let $A_1,\ldots,A_m\in\CP_4$ and
	$\omega\in\Delta_m$, and use \eqref{eq:CP-Y-data}. Set
	\begin{equation}\label{eq:CP-obstruction}
		\mathcal Q(\omega;\mathbf Y)
		=\sum_{k=1}^mw_k[Y_k,[M,Y_k]]\in\mathbb RE_{14}.
	\end{equation}
	Then the Karcher mean has the \emph{exact} formula
	\begin{equation}\label{eq:CP4-exact-Karcher}
	 \Lambda(\omega;\mathbf A)
			=\alpha\exp\left(M+\frac1{12}\mathcal Q(\omega;\mathbf Y)\right). 
	\end{equation}
	Furthermore,
	\begin{equation}\label{eq:CP4-exact-difference}
	 \Lambda(\omega;\mathbf A)-
			\mathrm{LE}(\omega;\mathbf A)
			=\frac{\alpha}{12}\mathcal Q(\omega;\mathbf Y). 
	\end{equation}
	In particular,
	\begin{equation}\label{eq:CP4-equality-iff}
		\Lambda(\omega;\mathbf A)=\mathrm{LE}(\omega;\mathbf A)
		\quad\Longleftrightarrow\quad
		\sum_kw_k[Y_k,[M,Y_k]]=0.
	\end{equation}
	For two inputs with weights $(w,1-w)$,
	\begin{equation}\label{eq:CP4-two-Q}
		\mathcal Q=w(1-w)[Y_2-Y_1,[Y_1,Y_2]].
	\end{equation}
\end{theorem}
\begin{proof}
	Scalar factors decouple from the Karcher equation, so we solve
	$\sum_kw_k\log(e^{-Z}e^{Y_k})=0$ for $Z\in\mathfrak{up}_4$.
	Applying \eqref{eq:CP-BCH3} gives
	\begin{equation}\label{eq:CP4-K-equation}
		\begin{split}
			0={}&M-Z-\frac12[Z,M]
			+\frac1{12}[Z,[Z,M]]
			+\frac1{12}\sum_kw_k[Y_k,[Z,Y_k]].
		\end{split}
	\end{equation}
	Modulo $\mathfrak n^3$ this equation yields $Z=M$:
	first $Z=M\pmod{\mathfrak n^2}$ and then
	$Z=M\pmod{\mathfrak n^3}$. Thus $Z=M+D$ for some
	$D\in\mathfrak n^3=\mathbb RE_{14}$.
	Since $D$ is central, every commutator containing $D$ vanishes.
	Equation \eqref{eq:CP4-K-equation} reduces to
	$0=-D+\mathcal Q/12$.
	This proves \eqref{eq:CP4-exact-Karcher}.
	Because $E_{14}$ annihilates every strictly upper triangular
	$4\times4$ matrix on both sides,
	$\exp(M+D)=\exp(M)+D$. This proves
	\eqref{eq:CP4-exact-difference} and the equality criterion.
	For $m=2$, expand $M=wY_1+(1-w)Y_2$ to obtain
	\eqref{eq:CP4-two-Q}.
\end{proof}

The following example gives a concrete nonzero difference on $\UP_4$. 

\begin{example}
	\label{ex:CP4-strict}
	Let $Y_1=E_{12}$, $Y_2=E_{23}+E_{34}$,
	$A_i=\exp Y_i$, and $w_1=w_2=1/2$.
	Then
	\[
	M=\tfrac12(E_{12}+E_{23}+E_{34}),\qquad
	\mathcal Q=-\tfrac14E_{14},
	\]
	and therefore
	\begin{equation}\label{eq:CP4-explicit-counterexample}
	 \Lambda(\tfrac12,\tfrac12;A_1,A_2)
			-\mathrm{LE}(\tfrac12,\tfrac12;A_1,A_2)
			=-\frac1{48}E_{14}\ne0. 
	\end{equation}
	This proves that dimension four is the first dimension in the
	family $\CP_n$ where the two means can differ. For every $n>4$,
embedding these inputs as $A_i\oplus I_{n-4}\in\UP_n$ gives the same
nonzero difference on the leading principal block.
\end{example}

\begin{corollary}[Cubic small-dispersion expansion]
	\label{cor:CPn-small-dispersion}
	Fix $Y_1,\ldots,Y_m\in\mathfrak{up}_n$ and let
	$U_k(\varepsilon)=e^{\varepsilon Y_k}$.
	Then, as $\varepsilon\to0$,
	\begin{equation*} %\label{eq:CPn-cubic-expansion}
		\begin{aligned}
			\log\Lambda(\omega;\mathbf U(\varepsilon))
			&=\varepsilon M+
			\frac{\varepsilon^3}{12}
			\sum_kw_k[Y_k,[M,Y_k]]+O(\varepsilon^4),\\
			\Lambda(\omega;\mathbf U(\varepsilon))
			-\mathrm{LE}(\omega;\mathbf U(\varepsilon))
			&=\frac{\varepsilon^3}{12}
			\sum_kw_k[Y_k,[M,Y_k]]+O(\varepsilon^4).
		\end{aligned}
	\end{equation*}
	The remainder is in any fixed matrix norm, with constants depending
	on the data and on $n$. For $n=4$ there is no remainder in the
	logarithmic formula, and the matrix difference is exactly cubic.
\end{corollary}
\begin{proof}
	By analytic dependence of the Karcher solution on the data,
	\[
	Z(\varepsilon)=\log\Lambda(\omega;
	e^{\varepsilon Y_1},\ldots,e^{\varepsilon Y_m})
	\]
	is analytic near zero.
	The degree-two truncation of BCH implies
	$Z(\varepsilon)=\varepsilon M+O(\varepsilon^3)$;
	indeed, the potential quadratic commutator term is
	$[M,M]=0$. Substitute this into the BCH expansion through
	homogeneous degree three. The Karcher equation then forces the
	cubic term to equal $\mathcal Q/12$.
	Exponentiating gives the second formula.
\end{proof}

\endgroup

\vspace{1cm}

\textbf{Acknowledgement}

There is no conflict of interest and all authors contributed equally.

%%%%%%%%%%%%%%%%%%%%%%%%%%%%%%%%%%%%%%%%%%%%%%%%%%

\end{document}